\documentclass[a4paper,12pt]{article}

\usepackage{amsmath,amsthm,amssymb,amscd}
\usepackage{pgfplots}
\usepackage{amssymb}
\usepackage{amsthm}
\usepackage{amsmath, amsfonts}
\usepackage{authblk}
\usepackage{epsfig}
\usepackage{amscd}
\usepackage{graphicx,wrapfig}
\usepackage{subcaption}
\usepackage[thinlines,thiklines]{easybmat}
\usepackage{hyperref}
\usepackage{capt-of}
\usepackage{float}
\usepackage{stackrel}
\usepackage{tikz-cd}
\def\natu           {\mathbb N}
\def\inte 		{\mathbb Z}

\def\D		{\mathbb D}
\def\F		{\mathbb F}
\def\A 		{{\cal A}}

\def\R		{{\cal R}}

\def\N		{{\cal N}}
\def\J		{{\cal J}}

\def\lra		{\longrightarrow}

\title{Lower central series of the Riordan group\\ over the field with two elements}

\author{Nikolai A. Krylov \thanks{Nikolai Krylov is supported by an AMS-Simons 
Research Enhancement Grant for Primarily Undergraduate Institution Faculty.}}
\affil{Department of Mathematics, Siena University, \newline 
515 Loudon Road, Loudonville NY 12211, USA \newline  
email: nkrylov@siena.edu}

\date{}

\begin{document}

\newtheorem{theorem}{Theorem}
\newtheorem{lemma}[theorem]{Lemma}
\newtheorem{claim}[theorem]{Claim}
\newtheorem{corollary}[theorem]{Corollary}
\newtheorem{conj}[theorem]{Conjecture}
\newtheorem{prop}[theorem]{Proposition}
\newtheorem{note}{Note}
\theoremstyle{definition}
\newtheorem{definition}[theorem]{Definition}
\newtheorem{example}[theorem]{Example}
\numberwithin{equation}{section}

\textheight=22cm
\topmargin=-4mm

\maketitle

\begin{abstract}
The Riordan group $\R$ over the field $\F_2$ is a split extension of the Appell 
subgroup by the Nottingham group $\N(\F_2)$. Using the lower central 
series of the Nottingham group obtained by C. Leedham-Green and S. McKay, 
the lower central series of $\R(\F_2)$ is calculated. Considering the Riordan 
group over an arbitrary commutative ring with identity, where all Riordan arrays 
have only 1s on the main diagonal, it is also proved that the abelianization of 
this group is isomorphic to the direct product of the abelianization of the 
corresponding Lagrange subgroup and the additive group of the ground ring.
\end{abstract}




\section{Introduction}

The Riordan group, introduced by L. Shapiro et al. \cite{Shapiro1}, is a group 
of infinite lower triangular matrices called Riordan arrays, whose columns 
consist of the coefficients of certain formal power series. For the detailed 
introduction to the subject, the reader is referred to the books by P. Barry 
\cite{Barry} and L. Shapiro et al. \cite{Shapiro2}, and a survey article by 
D. Davenport et al. \cite{Davenport}. Here are a few basic definitions 
and notations, which will be used in the paper. 

Let $\D[[t]]$ be the set all formal power series (f.p.s.) in indeterminate 
$t$ with coefficients in an arbitrary commutative ring $\D$ with identity, 
and denote its multiplicative group of units by $\D^*$. Take two f.p.s. 
$g(t) =\sum_{k=0}^\infty g_kt^k$ and $f(t) =\sum_{k=1}^\infty f_kt^k$, 
such that $g_0, f_1\in\D^*$, and for all $n,k\geq 0$, define
$$
d_{n,k} : = [t^n]g(t) f(t)^k,
$$ 
where $[t^n]h(t)$ stands for the coefficient of $t^n$ in the expansion of the 
f.p.s. $h(t)$. Then the Riordan array $(g,\,f)$ is defined to be the infinite 
lower triangular matrix
$$
A =(d_{n,k})_{n,k\geq 0}.
$$

Since each row of such a matrix has only finitely many nonzero terms, and 
$g_0,f_1\in\D^*$, these arrays form a group with the usual matrix multiplication, 
called the {\sl Riordan group over $\D$}. 
We will denote this group by $\R(\D)$, or simply by $\R$ 
when the ground ring $\D$ will be clear from the context. In terms of 
the pairs of f.p.s. the group operation is written as 
\begin{equation}
\label{Semdiprod}
\bigl(g_1(t) ,\,f_1(t)\bigr)\cdot \bigl(g_2(t) ,\,f_2(t)\bigr) = 
\bigl(g_1(t) g_2(f_1(t)),\,f_2(f_1(t))\bigr).
\end{equation}
The Riordan array $I = (1,\,t)$ is the group 
identity, and the inverse of the Riordan array $\bigl(g(t) ,\,f(t)\bigr)$ is the pair
$$
\bigl(g(t) ,\,f(t)\bigr)^{-1} = \left(\frac{1}{g(\bar{f}(t))},\,\bar{f}(t)\right),
$$
where we used the standard notation $\bar{f}(t)$ for the compositional 
inverse of $f(t)$, i.e. $\bar{f}(f(t)) = t$ and $f(\bar{f}(t)) = t$. 

As we see from \ref{Semdiprod}, $\R(\D)$ is 
a semidirect product of two proper subgroups: the Appell 
subgroup and the Lagrange (or associated) subgroup 
(\cite{Barry}, \cite{Shapiro2}). The Appell subgroup 
${\cal A}({\mathbb D})$ is abelian, normal, and consists of the 
Riordan arrays $(g(t) ,\,t)$. The Lagrange subgroup $\N({\mathbb D})$ 
is not commutative, and consists of the Riordan arrays $(1 ,\,f(t))$. 
When $\D^*$ has at least two elements, $\N(\D)$ properly contains 
the {\sl substitution group of formal power series} 
$\cal{J}({\mathbb D})$, which was studied in the 1950s by 
S. Jennings, \cite{Jennings} (see also review article by 
I. Babenko \cite{Babenko}). Elements of the substitution group, as 
Riordan arrays, are the pairs $(1 ,\,f(t))$, where the f.p.s. 
$f(t) =\sum_{k \geq 1} f_kt^k$ has the first coefficient $f_1 = 1$. The fact 
that the Riordan group is a semidirect product 
$\R(\D) \cong \A(\D) \ltimes \N(\D)$ 
can be written as a splitting short exact sequence 
\begin{equation}
\label{ES1}
\begin{tikzcd}
0 \arrow[r] & \A(\D) \arrow[r, "\mu"] & \R(\D) \arrow[r, "\rho"]  
& \N(\D) \arrow[l, bend left=20, "s"] \arrow[r] & 1,
\end{tikzcd}
\end{equation}
where $0$ indicates that $\A(\D)$ is commutative. Projection homomorphism 
$\rho$ is defined by $\rho\bigl((g(t),\, f(t))\bigr): =  (1,\, f(t))$, and 
the splitting homomorphism $s$ is the inclusion $s\bigl((1,\, f(t))\bigr):= (1 ,\, f(t))$.
These two maps clearly satisfy $\rho\circ s = id|_{\N(\D)}$. The 
homomorphism $\mu$ is an inclusion too.

When the ground ring $\D = \F_2 = \{0,1\}$, the Lagrange 
subgroup $\N(\F_2)$ is isomorphic to the Nottingham 
group, denoted by $\N(\F_2)$ as well. The lower central series 
of $\N(\F_2)$ was calculated by C.R. Leedham-Green 
and S. McKay in \cite{Leedham}, \S 12.4 (see also article 
\cite{Camina2} by R. Camina). Since I will use these 
calculations to obtain the lower central series of $\R(\F_2)$, 
I will follow closely their notations, which seem to be  
standardly accepted in the literature dedicated to the Nottingham 
group. Notice however that the corresponding notations for the Lagrange 
subgroup and truncated Lagrange subgroups used in the literature 
dedicated to the Riordan group are usually different. When 
$\D = \F_q$ and $q$ is odd, the Lagrange subgroup 
is also called the `full' group of f.p.s. under substitution, 
and contains the Nottingham group $\N(\F_q)$ as a proper 
normal subgroup, \cite{Camina2}. 

It was proved by A. Luz\'on et al. \cite{Luzon16}, that the Riordan group 
$\R(\D)$ is isomorphic to the inverse limit of the inverse system of Riordan 
matrices of finite size. Analogous result for the groups of formal power 
series under substitution was proved by D.L. Johnson in \cite{Johnson}.  
I will briefly recall the construction at the beginning 
of next section, since such Riordan matrices play a key role here. 
If the ring $\D$ is finite, the group $\R(\D)$ is isomorphic to the inverse limit of 
an inverse system of finite groups, and therefore is a {\sl profinite group}. A 
group $G$ is profinite if it is a compact Hausdorff topological group whose 
open subgroups form a base for the neighborhoods of the identity. Since 
for any finite ring $\D$, $\A(\D)$ is a closed normal subgroup of infinite index,
the Riordan group $\R(\D)$ is not {\sl just infinite}. Topology 
on the inverse limit is induced by the product topology on the infinite 
Cartesian product, where each finite group comes with 
the discrete topology. It is the weakest topology in which the projection 
map from the inverse limit group onto each finite group is continuous. 
For an introduction to the theory of profinite groups the reader is referred to 
one of the monographs \cite{Dixon}, \cite{Leedham}, \cite{Wilson}. 
When $\D = \F_{2}$, $\R(\D)$ is a pro-$2$ group, being the inverse 
limit of finite $2$-groups
$$
\R(\D) \cong\varprojlim\bigl(T\R_n(\D),~ P_{n-1}\bigr),
$$
where for each $n\geq 1$, the $2$-group $T\R_n(\D)$ consists of the Riordan arrays 
of finite size 
$$
\bigl(1 + \alpha_1t + \cdots + \alpha_nt^n,\, 
t + \beta_2t^2 + \cdots + \beta_nt^n \bigr), ~ \alpha_i,\beta_i \in \F_{2}, 1\leq i \leq n.
$$
The homomorphism $P_{n-1}: T\R_n\to T\R_{n-1}$ deletes the last row 
and column from the array $T\R_n$, i.e. 
$$
P_{n-1}\left(\bigl(1 + \alpha_1t + \cdots + \alpha_{n-1}t^{n-1} + \alpha_nt^n,\, 
t + \beta_2t^2 + \cdots + \beta_{n-1}t^{n-1} + \beta_nt^n \bigr)\right) 
$$
$$
= \bigl(1 + \alpha_1t + \cdots + \alpha_{n-1}t^{n-1},\, 
t + \beta_2t^2 + \cdots + \beta_{n-1}t^{n-1} \bigr).
$$

The main objective of this article is the lower central series of $\R(\F_2)$, i.e.  
$$
\R(\F_2) = \gamma_1(\R(\F_2)) \geq \gamma_2(\R(\F_2)) \geq \cdots,
$$
where $\gamma_i(G) = [\gamma_{i-1}(G),\,G]$ for all $i\geq 2$. Our result is 
stated in Theorem \ref{LCSRior2} in terms of the lower central series of the  
Nottingham group $\N(\F_2)$, which is given in \cite{Leedham}, \S 12.4. 
In particular, $\N^{ab}(\F_2)\cong \inte_2\times \inte_2 \times \inte_4.$
When $p$ is odd, the corresponding computations 
for $\N(\F_p)$ are given in Theorem 2, \cite{Camina2}.
It is important to notice that by a subgroup of a topological group $G$ one means 
a nonempty closed subset of the topological space $G$, which is also closed under 
multiplication and inversion. The commutator of two closed subgroups of a 
topological group is not automatically closed as the topological subset, 
therefore to have a well defined lower central series, one needs the notion 
of a {\sl topologically generating set}. In particular, if $X$ and $Y$ are normal 
closed subgroups of $\R(\D)$, we define $[X,\,Y]$ as the smallest closed 
subgroup containing all the commutators $[x,\,y] = x^{-1}y^{-1}xy$, where 
$x\in X$ and $y\in Y$ (see \cite{Jennings}, \S 2). That  is how we will understand 
all commutators $[X,\,Y]$ here, unless it will be specifically stated otherwise. 

Using results of R. Camina \cite{Camina2}, and C.R. Leedham-Green 
and S. McKay \cite{Leedham}, on the lower central 
series of the Nottingham group over a finite field with characteristic 
$p> 2$, Gi-S. Cheon, N.-P. Chung, and M.-N. Phung calculated in 
\cite{Cheon1} the lower central series of the corresponding Riordan 
group. S.I. Bogataya and S.A. Bogatyi utilizing extensive calculations 
in \cite{Bogataya} found the high commutants and gave explicit 
abelianization homomorphism of the Jennings - Lagrange - Nottingham 
group $\N(\F_2)$. D.L. Johnson, \cite{Johnson} calculated the lower 
central series of the group of f.p.s. under substitution over the rationals, 
and I.K. Babenko and S.A. Bogatyi in \cite{BabenkoBogat}
found the abelianization of the group of f.p.s. under substitution  
with integer coefficients. Recently, A. Luz\'on, M.A. Mor\'on, 
and L.F. Prieto-Mart\'inez calculated the derived series of the Riordan 
group over a field of characteristic 0, \cite{Luzon23}. The present work 
is motivated by this research and offers further insights into the algebraic 
structure of the Riordan group. 

The paper is organized as follows. In section 2, we will review the inverse limit 
description of the Riordan group, and use the short exact sequence \ref{ES1} 
to derive a commutative diagram connecting the groups $S\R$ and $S\N$ 
with their abelianizations. Symbol ``$S$" indicates that all the corresponding 
matrices have 1s on the main diagonal. From this diagram  
we obtain an isomorphism $S\R^{ab}(\D) \cong \D\times S\N^{ab}(\D)$,  
for an arbitrary commutative ring $\D$ with identity. This isomorphism, together 
with calculations of I. Babenko and S. Bogatyi, \cite{BabenkoBogat}, gives the 
answer to a question posed by L. Shapiro in \cite{Shapiro3} about the commutator 
group of the subgroup of $\R(\inte)$, 
that has all 1s on the main diagonal (see Corollary \ref{ShapiroAnswer} below). 
In section 3 we will classify all truncated Appell subgroups over $\F_2$ as finite 
abelian groups, and prove our main result, Theorem \ref{LCSRior2}. 
When $\D = \F_2$, $\gamma_n(\R) = \A_{2n-3} \ltimes \gamma_n(\N)$, and
$$
\gamma_i(\R)/\gamma_{i+1}(\R) \cong
\left\{ \begin{array}{cl}
(\inte_2)^3\times \inte_4 & \mbox{if $i = 1$ } \\
(\inte_2)^4 & \mbox{if $i > 1$ is even } \\
(\inte_2)^6 & \mbox{if $i > 1$ is odd. } \\
\end{array}
\right.
$$
At the end, we will show how to embed the dihedral group $D_{2^{n+1}},\forall n\geq 0$ 
into the truncated Riordan group $T\R_{2^n}(\F_2)$.


\section{Short exact sequences}

Let $\D$ be an arbitrary, fixed, commutative ring with identity 1. 
For any $n\geq 0$, take the 
general linear group of all invertible $(n+1)\times (n+1)$ matrices over 
$\D$, ${\rm GL}(n+1,\D)$, and denote by ${T\R_n(\D)} : = \Pi_n(\R(\D))$ 
the image of a natural homomorphism (called {\sl truncation})
\begin{equation}
\label{Pin}
\Pi_n: \R(\D) \longrightarrow {\rm GL}(n+1,\D)
\end{equation}
defined by
$$
\Pi_n\bigl((d_{i,j})_{i,j\geq 0}\bigr) = (d_{i,j})_{0\leq i,j\leq n}.
$$
Deleting the last row and column from $T\R_{n+1}(\D)$ produces 
$T\R_n(\D)$, and this way we obtain another natural homomorphism, 
which will be denoted by $P_n: T\R_{n+1}(\D) \longrightarrow T\R_n(\D)$. 
In terms of the formal notations we have 
$$
P_n\bigl((d_{i,j})_{0\leq i,j\leq n+1}\bigr) = (d_{i,j})_{0\leq i,j\leq n},
$$
and clearly $\Pi_n = P_n\circ\Pi_{n+1}$. The homomorphism 
$P_n$ is onto, and for all $n\geq 0$, the homomorphisms 
$\Pi_{n+1}, \Pi_n$ and $P_n$ make the commutative triangle.
\begin{equation}
\label{ComD}
\begin{tikzcd}
 & \R(\D) \arrow[dl,"\Pi_{n+1~}" left, near start] \arrow[dr, "\Pi_n"] & \\
T\R_{n+1}(\D) \arrow{rr}{P_n} & & T\R_n(\D)\\
\end{tikzcd}
\end{equation}
Finite size Riordan arrays $T\R_n(\D)$ can be constructed iteratively using 
the $A$ and $Z$ sequences, and it was shown in \cite{Luzon16}, 
that the Riordan group $\R(\D)$ is isomorphic to the inverse limit 
$$
\R(\D) \cong\varprojlim\bigl(T\R_{n+1}(\D),\,P_n\bigr)
$$
(c.f. Proposition 1. in \cite{Johnson}, for the group of f.p.s. under substitution).

Since $\R(\D) \cong \A(D)\ltimes \N(\D)$, restricting the homomorphism 
$\Pi_n$ onto the subgroups $\A(\D)$ and $\N(\D)$ will produce 
truncated Appell and Lagrange subgroups, and for every $n\geq 0$, 
we get $T\R_n(\D) \cong T\A_n(D)\ltimes T\N_n(\D)$. 
Truncated group $T\N_n(\D)$ was also introduced by R. Camina in 
\cite{Camina2}, as the factor group of $\N(\D)$ by the kernel of 
$$
\Phi_n: \bigl(f : t\mapsto \sum\limits_{i =1}^{\infty} \alpha_it^i\bigr)
\mapsto \bigl(f_n : t\mapsto \sum\limits_{i =1}^n \alpha_it^i\bigr).
$$ 
Clearly, the maps $\Phi_n$ and  
$\Pi_n\left|_{\N(\D)} \right.$ are the same after identifying the 
f.p.s. $\sum\limits_{i =1}^{\infty} \alpha_it^i$ with the Riordan 
array $\bigl(1,\,\sum\limits_{i =1}^{\infty} \alpha_it^i\bigr)$ . The kernel of 
$\Phi_n$ is a normal subgroup of $\N(\D)$ denoted by $\N_n(\D)$, so 
$T\N_n(\D) \cong \N(\D)/\N_n(\D)$. When $\D$ is finite, these normal 
subgroups $T\N_n$ form a countable basis for the neighborhoods of 
the identity in $\N(\D)$. Analogously, when $\D$ is finite, 
the kernels of $\Pi_n(\R(\D)),n\geq 0$ form a countable basis for the 
neighborhoods of the identity $(1,\,t)\in \R(\D)$. 

It will be convenient to use certain 
elements of the group $\N_n(\D)$, so let us recall some notations 
from \cite{Camina2} and \cite{Leedham}. The Riordan array 
$(1,\, t + \alpha t^{n+1})$ will be denoted by $e_n[\alpha]$, or simply 
by $e_n$ if $\alpha=1$, and $e_0 = (1,\,t)$ will stand for the identity 
element. Sequence \ref{ES1} naturally induces the following split 
short exact sequence of the truncated groups for all $n\in\natu$ 
\begin{equation}
\label{ES2}
\begin{tikzcd}
0 \arrow[r] & T\A_n(\D) \arrow[r, "\mu_n"] & T\R_n(\D) \arrow[r, "\rho_n"]  
& T\N_n(\D) \arrow[l, bend left=20, "s_n"] \arrow[r] & 1,
\end{tikzcd}
\end{equation}
where maps $\mu_n$, $\rho_n$, and $s_n$ are the restrictions 
of their corresponding analogs from \ref{ES1}. 
If we have a f.p.s. $g(t)$, we can take $g_n(t)$, the $n$-th degree 
Taylor polynomial of it for any $n\geq 0$. Truncated Riordan arrays 
$\bigl(g_n(t),\, t\bigr)$ are also known as Toeplitz matrices of size 
$(n + 1) \times (n + 1)$ with zeroes above the main diagonal. 
Restriction of the homomorphism $P_n$ onto the Appell 
subgroup $T\A_{n+1}(\D)$ gives an epimorphism
$$
P_n\left |_{T\A_{n+1}(\D)}\right.: T\A_{n+1}(\D) \longrightarrow T\A_n(\D),
$$
with $T\A_0(\D) \cong \D^*$. Since a Toeplitz matrix with zeroes above 
the main diagonal is completely determined by its most left column, 
$\ker(P_n\left |_{T\A_{n+1}(\D)}\right.)\cong \D$ for all $n\geq 0$. 
Hence, there is the exact sequence 
\begin{equation}
\label{ES22}
\begin{CD}
0 @>>> \D @>>> {T\A_{n+1}(\D)} @>{P_n}>> {T\A_n(\D)} @>>> 0,
\end{CD}
\end{equation}
\\
where for brevity, $P_n$ was also used to denote the restriction 
homomorphism $P_n\left |_{T\A_{n+1}(\D)}\right.$. Notice that $\D$ 
here (and later in a similar context) stands for the additive group of the 
ring, i.e. $\D^+$. We also have $T\A_1(\D)\cong \D\times \D^*$, where
the isomorphism is given by the map
$$
\begin{pmatrix}
u & 0\\
x & u
\end{pmatrix} \longmapsto (xu^{-1},\, u)\in \D\times \D^*.
$$
In general however, $T\A_n(\D)\ncong \D^n\times \D^*$. Since the matrix 
$$
\begin{pmatrix}
1 & 0 & 0\\
1 & 1 & 0\\
0 & 1 & 1
\end{pmatrix}
$$
corresponding to $(1 + t,\,t)$, has order 4 in $T\A_2(\F_2)$, we have 
 ${T\A_2(\F_2)}\cong C_4$, the cyclic group of order 4. We will discuss 
the group structure of $T\A_n(\F_2)$ in details in the next section. 
Similarly, restricting $P_n$ onto the subgroup $T\N_{n+1}(\D)$, one 
obtains the short exact sequence 
\begin{equation}
\label{ES3}
\begin{CD}
0 @>>> \D @>>> {T\N_{n+1}(\D)} @>{P_n}>> {T\N_n(\D)} @>>> 1,
\end{CD}
\end{equation}
with $T\N_1(\D) \cong \D^*$ and $T\N_2(\D) \cong  \D\ltimes \D^*$. 


When $n\geq 1$, it is convenient to connect together three exact sequences 
\ref{ES2}, \ref{ES22}, and \ref{ES3}, in one commutative diagram (see Lemma 1. 
and diagram (2.6) in the article \cite{Krylov2} by T.-X. He and the author).

\begin{equation}
\label{BCD1}
\begin{tikzcd}[column sep=4ex, row sep=4ex]
& 0 \arrow[d] & 0 \arrow[d] & 0 \arrow[d] & \\
0 \arrow[r] & \D \arrow[r] \arrow[d] & \D\times \D \arrow[r] \arrow[d] & 
\D \arrow[l, bend left=20] \arrow[r] \arrow[d] & 0 \\
0 \arrow[r] & T\A_{n+1}(\D) \arrow[r, "\mu_{n+1}"] \arrow[d] 
& T\R_{n+1}(\D) \arrow[r, "\rho_{n+1}"]  \arrow[d] 
& T\N_{n+1}(\D) \arrow[l, bend left=20, "s_{n+1}"] \arrow[r] \arrow[d] & 1 \\
0 \arrow[r] & T\A_n(\D) \arrow[r, "\mu_n"] \arrow[d] & T\R_n(\D) \arrow[r, "\rho_n"]  \arrow[d]
& T\N_n(\D) \arrow[l, bend left=20, "s_n"] \arrow[r] \arrow[d] & 1 \\
& 0 & 1 & 1 & \\
\end{tikzcd}
\end{equation}

L. Shapiro asked about the commutator group of the subgroup of 
$\R(\inte)$, that has all 1s on the main diagonal (see \cite{Shapiro3}, Q10). 
He denoted such a subgroup by $S\R$. Using results of I. Babenko and 
S. Bogatyi \cite{BabenkoBogat}, on the commutator subgroup 
$[{\cal J}(\inte),\,{\cal J}(\inte)]$, of the group of f.p.s. under substitution 
${\cal J}(\inte)$, we answer Shapiro's question below in Corollary 
\ref{ShapiroAnswer}. Motivated by this problem, instead of 
the Lagrange subgroup $\N(\D)$ and its truncated groups $T\N_n(\D)$, 
from now on I will consider the group of f.p.s. under substitution 
with coefficients in $\D$, that is ${\cal J}(\D)$, and the corresponding 
truncated groups $T{\cal J}_n(\D)$. Exact sequence \ref{ES3} holds true for 
such groups, but $T{\cal J}_1(\D) = \{e\}$, and $T{\cal J}_2(\D) \cong \D$, 
since we have all 1s on the main diagonal of every array from ${\cal J}(\D)$. 
Similarly, all formal power series $g(t)$ used in the S-Appell subgroup 
will have $g(0) = 1$. Thus, the main diagonals of all the corresponding 
Riordan arrays will have all 1s, and to follow Shapiro's notation, I will 
add ``S" to all the notations for such groups, and their truncated versions, 
i.e. $S\A,~TS\A_n$ and $S\R,~TS\R_n$. For brevity, I will also drop 
the symbol $\D$ till Theorem \ref{Isom1}, and subsection 2.1. 
Furthermore, to be in agreement 
with the notation $e_n[\alpha] = (1,\,t+\alpha t^{n+1})$, I will use 
$a_k(\beta)$ for the Riordan array $(1 + \beta t^k,\,t)$, where 
$\beta\in \D$. In such notations the kernel of the projection 
$P_n : TS\R_{n+1}\to TS\R_n$ is generated by $a_{n+1}(\beta)$ and 
$e_{n}[\alpha]$, for some elements $\alpha, \beta$ from the ground ring. 

Recall that for the groups consisting of 
finite size Riordan arrays, the discrete topology is used, and  
the commutator subgroup $[G,\,G]$ is considered as the 
{\sl abstract} subgroup of $G$. In particular, every element of 
$[G,\,G]$ is a product of finitely many commutators.

We will show next that 
for all $n\geq 1$, the abelianization of the truncated Riordan group 
$TS\R_n^{ab}$ is isomorphic to the direct product of $\D$ 
with the abelianization of the truncated subgroup 
$T{\cal J}_n^{ab}$ (Theorem \ref{Isom1}). Next lemma  
proves that for any $\alpha\in\D$, and any integer $k\geq 2$, the 
array $a_k(\alpha)$ can be written as a commutator, i.e. 
$(1 + \alpha t^k, \, t) = [x,\,y]$, and hence 
$a_k(\alpha)\in \gamma_2(S\R)$. 

\begin{lemma}
\label{lemCOM1}
For any $\alpha\in\D$ and $k\geq 2$, $a_k(\alpha) \in \gamma_2(S\R)$. 
We also have $a_1(\alpha) \notin \gamma_2(TS\R_n)$ for all $n\geq 1$, unless 
$\alpha = 0$.
\end{lemma}
\begin{proof}
Since $\D$ is commutative, the composition of homomorphisms
$$
P_2\circ \cdots \circ P_{n} : TS\R_{n+1} \to TS\R_1 = 
\{a_1(\alpha)~|~\alpha\in \D\} 
$$
maps any commutator to the identity,
so $a_1(\alpha)\notin \gamma_2(TS\R_n)$, unless it is the identity. 
Replacing the composition by $\Pi_1$, one shows 
that $a_1(\alpha)\notin \gamma_2(S\R)$.
To prove the first statement, take polynomials $g(t) = 1 + t$, and 
$f(t) = t + \alpha t^k + \alpha t^{k+1}, k\geq 2$. 
Then direct computations show 
$$
\bigl[\bigl(g(t),\,t\bigr),\, \bigl(1,\,\bar{ f}(t)\bigr)\bigr] = 
\left(\frac{1}{1 + t},\, f(t)\right) \left(1 + t,\,\bar{f}(t)\right) 
$$
$$
= \left(\frac{1 + t + \alpha t^k + \alpha t^{k+1}}{1 + t},\,t\right) 
= \bigl(1 + \alpha t^k,\, t\bigr),
$$
as required.
\end{proof}

\noindent If we define subgroups $\A_n$ in a similar way to the 
subgroups $\N_n$, i.e. as 
\begin{equation}
\label{NormAppell}
\A_n:= \ker\left(\left.\Pi_n\right|_{\A}: \A\lra T\A_n\right),
\end{equation}
Lemma \ref{lemCOM1} immediately implies that 
$\A_1 = \A\cap \gamma_2(S\R) = S\A\cap \gamma_2(S\R)$.
Notice that the main diagonals of elements from $\A_n$ consist of all 1s, and 
$\A_n = S\A_n$ for all $n\geq 0$, where $\A_0 = S\A \lneq \A$. Therefore, 
considering $\gamma_2(S\R)$ as the discrete subgroup of $S\R$, we obtain  
\begin{corollary}
\label{corCOM2}
For every integer $n\geq 1$, there is an isomorphism
$$
TS\A_n /\bigl(TS\A_n\cap \gamma_2(TS\R_n)\bigr)
\cong \left\{
\begin{pmatrix}
1 & 0\\
\alpha & 1 \\ 
\end{pmatrix} ~ |~  \alpha\in \D\right\} = TS\A_1 \cong \D,
$$
and also 
$$
S\A /\bigl(S\A\cap \gamma_2(S\R)\bigr) = S\A/\A_1\cong TS\A_1 \cong \D.
$$
\end{corollary}

If one applies the abelianization functor to the short exact sequence 
\begin{equation}
\label{ESS1}
\begin{tikzcd}
0 \arrow[r] & TS\A_n \arrow[r, "\mu_n"] & TS\R_n \arrow[r, "\rho_n"]  
& T{\cal J}_n \arrow[l, bend left=20, "s_n"] \arrow[r] & 1,~~ n\geq 1
\end{tikzcd}
\end{equation}
one derives the following commutative diagram,
where the projections $F_{\R}$ and $F_{\J}$ are the abelianization 
homomorphisms. Recall that the abelianization functor is right exact, but 
it is straightforward to check the exactness at each group of \ref{BCD2} 
by standard diagram chasing. Explanation of all maps is given below.

\begin{equation}
\label{BCD2}
\begin{tikzcd}
&  & 1 \arrow[d] & 1 \arrow[d] & \\
&  & 
\gamma_2(TS\R_n) \arrow[r, "\rho^{(2)}_n"] \arrow[d, "\iota_{S\R}"] & 
\gamma_2(T\J_n) \arrow[l, bend left=20, "s^{(2)}_n"] \arrow[d, "\iota_{\J}"] \arrow[r] & 1 \\
0 \arrow[r] & TS\A_n \arrow[r, "\mu_n"] \arrow[d, "F_{S\A}"] 
& TS\R_n \arrow[r, "\rho_n"]  \arrow[d, "F_{S\R}"] 
& T\J_n \arrow[l, bend left=20, "s_n"] \arrow[r] \arrow[d, "F_{\J}"] & 1 \\
0 \arrow[r] & \ker(\rho^{ab}_n) \arrow[r, "\mu^{ab}_n"]  & 
TS\R_n^{ab} \arrow[r,"\rho^{ab}_n"]  \arrow[d] & 
T\J_n^{ab} \arrow[l, bend left=20, "s^{ab}_n"] \arrow[r] \arrow[d] & 0 \\
&  & 0 & 0 & \\
\end{tikzcd}
\end{equation}

\noindent Take the epimorphism $\rho_n$ and restrict it onto the commutator 
subgroup to get the epimorphism $\rho_n^{(2)}$. Then $\rho_n^{ab}$ is the 
induced homomorphism between the abelianizations, and since $F_{\J}$ 
and $\rho_n$ are epimorphisms, $\rho_n^{ab}$ is onto as well. In a similar way, 
take the splitting monomorphism $s_n$ and restrict it onto $\gamma_2(T\J_n)$ to 
obtain the monomorphism $s_n^{(2)}$, which satisfies 
$$
\rho^{(2)}_n\circ s^{(2)}_n = Id_{\gamma_2(T\J_n)}.
$$ 
The homomorphism $s_n^{ab}: T\J_n^{ab} \to TS\R_n^{ab}$ is the induced 
homomorphism between the quotient groups. Commutativity 
$s_n^{ab}\circ F_{\J} = F_{S\R}\circ s_n$ follows directly from the 
definition of $s_n^{ab}$. The fact that $s_n^{ab}$ satisfies  
$\rho_n^{ab}\circ s_n^{ab} = Id_{T\J_n^{ab}}$, follows from the commutativity 
$\rho_n^{ab} \circ F_{S\R} = F_{\J} \circ \rho_n$. Indeed, take any 
$\bar{x}\in T\J_n^{ab}$ and any $x\in T\J_n$ such that $F_{\J}(x) = \bar{x}$. 
Then according to the definition of $s_n^{ab}$, 
$$
\rho_n^{ab}\circ s_n^{ab}(\bar{x}) = \rho_n^{ab} \circ F_{S\R}\circ s_n(x)
= F_{\J} \circ \rho_n \circ s_n(x) = F_{\J}(x) = \bar{x}.
$$
In particular, the splitting map $s_n^{ab}$ is a monomorphism too.

Monomorphisms $\iota_{S\R}$ and $\iota_{\J}$ are the inclusion homomorphism 
corresponding to the abelianization projections $F_{S\R}$ and $F_{\J}$, 
and the map $F_{S\A}$ is defined as the composition 
$F_{S\A}:= F_{S\R}\circ \mu_n$, which is the same as the restriction 
of $F_{S\R}$ onto the subgroup $TS\A_n$. 
To see why $F_{S\A}$ is onto, one can use 
diagram chasing as follows. Take any $x\in \ker(\rho^{ab}_n)$, and any 
$\bar{x}\in TS\R_n$ s.t. $F_{S\R}(\bar{x}) = x$. 
Since $F_{\J}\circ \rho_n(\bar{x}) = 0$,  $\exists y\in \gamma_2(TS\R_n)$ s.t. 
$\iota_{\J}\circ \rho^{(2)}_n(y) = \rho_n(\bar{x})$. Then 
$\rho_n\bigl(\iota_{S\R}(y^{-1})\cdot \bar{x}\bigr) = 1$, i.e. 
$\iota_{S\R}(y^{-1})\cdot \bar{x}\in TS\A_n$, and 
$F_{S\A}\bigl(\iota_{S\R}(y^{-1})\cdot \bar{x}\bigr) = 
F_{S\R} \circ \mu_n\bigl(\iota_{S\R}(y^{-1})\cdot \bar{x}\bigr) = 
F_{S\R} (\bar{x}) = x$. Map $\mu^{ab}_n$ is the 
inclusion $\ker(\rho^{ab}_n)\hookrightarrow TS\R_n^{ab}$. 
The kernel of $F_{S\A}$ equals 
the intersection of $TS\A_n$ with the kernel of $F_{S\R}$, i.e. 
$$
\ker(F_{S\A}) = TS\A_n\cap \gamma_2(TS\R_n).
$$
We saw in Corollary \ref{corCOM2} that 
$TS\A_n/(TS\A_n\cap \gamma_2(TS\R_n)) \cong \D$, 
so $\ker(\rho_n^{ab})\cong \D$. It is clear that 
$TS\A_n\cap \gamma_2(TS\R_n)\leq \ker(\rho^{(2)}_n)$. If we take any 
$x\in \ker(\rho^{(2)}_n) < \gamma_2(TS\R_n)$, then 
$\rho_n\circ \iota_{S\R}(x) = 0$, 
but $\iota_{S\R}$ is an inclusion, so we must have 
$$
x  \in \ker (\rho_n) = TS\A_n,
$$
i.e. $x\in TS\A_n\cap \gamma_2(TS\R_n)$, and thus 
$TS\A_n\cap \gamma_2(TS\R_n) = \ker(\rho^{(2)}_n)$.

If we denote this intersection $TS\A_n\cap \gamma_2(TS\R_n)$ by 
$\widetilde{TS\A}_n$, then we can put all the information from corollary 
\ref{corCOM2} and the diagram \ref{BCD2} into 
the following commutative diagram, where $n\geq 1$.

\begin{equation}
\label{BCD3}
\begin{tikzcd}
& 0 \arrow[d] & 1 \arrow[d] & 1 \arrow[d] & \\
0 \arrow[r] & \widetilde{TS\A}_n \arrow[r, "\mu^{(2)}_n"] \arrow[d] & 
\gamma_2(TS\R_n) \arrow[r, "\rho^{(2)}_n"] \arrow[d, "\iota_{S\R}"] & 
\gamma_2(T\J_n) \arrow[l, bend left=20, "s^{(2)}_n"] \arrow[d, "\iota_{\J}"] \arrow[r] & 1 \\
0 \arrow[r] & TS\A_n \arrow[r, "\mu_n"] \arrow[d, "F_{S\A}"] 
& TS\R_n \arrow[r, "\rho_n"]  \arrow[d, "F_{S\R}"] 
& T\J_n \arrow[l, bend left=20, "s_n"] \arrow[r] \arrow[d, "F_{\J}"] & 1 \\
0 \arrow[r] & \D \arrow[r, "\mu^{ab}_n"] \arrow[d] & 
TS\R_n^{ab} \arrow[r,"\rho^{(ab)}_n"]  \arrow[d]
& T\J_n^{ab} \arrow[l, bend left=20, "s^{ab}_n"] \arrow[r] \arrow[d] & 0 \\
& 0 & 0 & 0 & \\
\end{tikzcd}
\end{equation}

\noindent When a short exact sequence of abelian groups 
$$
\begin{tikzcd}
0 \arrow[r] & H \arrow[r] & G \arrow[r] & Q \arrow[r] \arrow[l, bend left=30] & 0
\end{tikzcd}
$$
splits, $G\cong H\times Q$. Therefore the 
bottom exact sequence in \ref{BCD3} implies that for all $n\geq 1$, 
the abelianization of the truncated S-Riordan group is isomorphic 
to the direct product of the abelianization of the truncated 
subgroup $\J_n$ with $\D$, where the elements of $\D$ 
correspond to the elements of $TS\A_1(\D) = TS\R_1(\D)$. 

Moreover, if we start with the exact sequence \ref{ES1} instead of 
\ref{ES2}, and follow the same approach we used to derive \ref{BCD3}, 
we obtain a completely analogous commutative diagram 
\ref{BCD3ORIG}, where all truncated groups are replaced by 
the corresponding infinite groups $S\A,~S\R$, and $\J$
considered as abstract groups (i.e. with the discrete topology). 
All the maps in \ref{BCD3ORIG} are defined identically to the 
corresponding maps in \ref{BCD3}, and by Corollary \ref{corCOM2},
$S\A\cap \gamma_2(S\R) = \A_1$.

\begin{equation}
\label{BCD3ORIG}
\begin{tikzcd}
& 0 \arrow[d] & 1 \arrow[d] & 1 \arrow[d] & \\
0 \arrow[r] & \A_1 \arrow[r, "\mu^{(2)}"] \arrow[d] & 
\gamma_2(S\R) \arrow[r, "\rho^{(2)}"] \arrow[d, "\iota_{S\R}"] & 
\gamma_2(\J) \arrow[l, bend left=20, "s^{(2)}"] \arrow[d, "\iota_{\J}"] \arrow[r] & 1 \\
0 \arrow[r] & S\A \arrow[r, "\mu"] \arrow[d, "F_{S\A}"] 
&S\R \arrow[r, "\rho"]  \arrow[d, "F_{S\R}"] 
& \J \arrow[l, bend left=20, "s"] \arrow[r] \arrow[d, "F_{\J}"] & 1 \\
0 \arrow[r] & \D \arrow[r, "\mu^{ab}"] \arrow[d] & 
S\R^{ab} \arrow[r,"\rho^{(ab)}"]  \arrow[d]
& \J^{ab} \arrow[l, bend left=20, "s^{ab}"] \arrow[r] \arrow[d] & 0 \\
& 0 & 0 & 0 & \\
\end{tikzcd}
\end{equation}

From the bottom exact sequence in \ref{BCD3ORIG}, we obtain the
abelianization of the S-Riordan group in terms of the abelianization 
$\J^{ab}(\D)$. Let us summarize this discussion as

\begin{theorem} 
\label{Isom1}
Let $\D$ be an arbitrary commutative ring with identity.
\begin{itemize}
\item For each $n\in\natu$,
$$
TS\R_n^{ab}(\D)\cong TS\A_1(\D)  \times T\J_n^{ab}(\D) 
\cong \D \times T\J_n^{ab}(\D).
$$
\item Considering $S\R(\D)$ and $\J(\D)$ as abstract groups,
$$
S\R^{ab}(\D) \cong \D \times \J^{ab}(\D).
$$
\end{itemize}
\end{theorem}

\subsection{Abelianization of $S\R(\inte)$}

Now we are in position to answer the question, which L. Shapiro 
asked in \cite{Shapiro3}. 
For the topology on ${\cal J}(\inte)$ the reader is referred 
to \cite{BabenkoBogat}, or \S 5.8 in Babenko's survey 
\cite{Babenko}, and topology on $S\R(\inte)$ is naturally induced 
by the epimorphism $S\R(\inte) \to {\cal J}(\inte)$. Here is the 
statement of Babenko-Bogatyy's theorem (see \cite{BabenkoBogat}, 
Proposition 2.9, or \cite{Bogataya}, Theorem 1.6).

~

\noindent {\bf Proposition 2.9.} The following conditions on 
$$
f = t\bigl(1 + \alpha_1t + \alpha_2t^2+\alpha_3t^3 + \alpha_4t^4+\alpha_5t^5 
+ \alpha_6t^6 + {\cal O}(t^7) \bigr)\in {\cal J}(\inte)
$$
are equivalent:
\begin{itemize}
\item[]{1)} $\alpha_1 = \alpha_2 = 0 ~~ \mbox{and} ~~ 
\alpha_3\equiv\alpha_4\equiv\alpha_6 \pmod{2},$
\item[]{2)} $f$ belongs to the kernel $\ker(\varphi)$ of a surjective homomorphism 
$\varphi: {\cal J}(\inte)\to \inte\oplus \inte\oplus\inte_2\oplus\inte_2$, 
which is defined in Proposition 2.1, \cite{BabenkoBogat}, 
\item[]{3)} $f$ belongs to $\bigl[{\cal J}(\inte),\,{\cal J}(\inte)\bigr]$, that is, $f$ is 
the limit of a finite product of commutators of elements of ${\cal J}(\inte)$.
\end{itemize}

\noindent If we Identify the f.p.s. $f(t)$ with the Riordan array $(1,\,f(t))$, this 
Proposition 2.9. together with Theorem \ref{Isom1} imply  

\begin{corollary}
\label{ShapiroAnswer}
The Riordan array with integer coefficients $c_i,\alpha_j\in\inte$ 
$$
(g,\,f) = \bigl(1 + \sum\limits_{i\geq 1} c_i t^i,\, 
t(1 + \alpha_1 t + \cdots + \alpha_6t^6 +{\cal O}(t^7))\bigr) 
\in \bigl[S\R(\inte),\,S\R(\inte)\bigr], 
$$
if and only if $c_1 = 0$, and $\alpha_1 = \alpha_2 =0$ and 
$\alpha_3\equiv \alpha_4\equiv \alpha_6 \pmod{2}$.\\ In particular,
$$
S\R^{ab}(\inte) \cong (\inte)^3\times (\inte_2)^2.
$$
\end{corollary}


\section{When $\D=\F_2$}

In this section we describe the lower central series 
$\gamma_i(\R(\F_2)),\,i\geq 1$ together with the lower central 
quotients of $\R(\F_2)$, in terms of the lower central 
series and the corresponding quotients of the Nottingham group 
$\N(\F_2)$. The latter were calculated in \cite{Leedham},  
Proposition 12.4.30 (c.f. \cite{Camina2}, Theorem 14.). 
S. Bogataya and S. Bogatyi also obtained the commutator 
$\bigl[\N(\F_2),\,\N(\F_2)\bigr]$ together with the abelianization 
homomorphism in \cite{Bogataya}.
Notice that $\J(\F_2) \equiv \N(\F_2)$, and 
$SG(\F_2) \equiv G(\F_2)$, where $G\in\{\A, \A_n, \R,\R_n\}$. 
For the rest of the paper all Riordan arrays will have coefficients 
in $\F_2$, so I will drop the ground field $\D = \F_2$ from all the 
corresponding notations. Recall also that the group $\A_n$ 
consists of the infinite Riordan arrays 
$$
\bigl(1 + \alpha_{n+1}t^{n+1} +  \alpha_{n+2}t^{n+2}  + \cdots, \,t\bigr),
$$ 
and in particular, $\A_0 \equiv S\A \equiv \A$. Truncated Appell subgroups $T\A_n$
are finite abelian 2-groups, and we can classify them in terms of the 
decomposition 
\begin{equation}
\label{Aclass}
T\A_n \cong \inte_{2^{r_1}} \times \inte_{2^{r_2}} \times \ldots \times \inte_{2^{r_k}},
\end{equation}
where $2^{r_1},~2^{r_2},~ \ldots ,~2^{r_k}$ are the {\sl invariant factors} of $T\A_n$, i.e. 
$r_i \geq  r_{i + 1}$ for each $i\geq 1$. It is clear that the cardinality $|T\A_n| = 2^n$. 
As we saw above, the sequence 
\begin{equation}
\label{AESZ2}
\begin{CD}
0 @>>> \inte_2 @>>> {T\A_{n+1}(\F_2)} @>{P_n}>> {T\A_n(\F_2)} @>>> 0
\end{CD}
\end{equation}
doesn't split in general, since $T\A_2(\F_2) \cong C_4$ is a cyclic group 
of order 4. It will follow from the proof below that \ref{AESZ2} splits 
if and only if $2~|~n$.

\begin{theorem}
\label{DecompA}
For each integer $n\geq 1$,
\begin{equation}
\label{DecompA2}
T\A_{n+1} \cong \inte_{2^{r_1}} \times \inte_{2^{r_2}} \times \ldots \times \inte_{2^{r_k}}, 
\end{equation}
where
$$
k = \left\lfloor{\frac{n+1}{2}}\right\rfloor, ~ \mbox{and} ~ 
r_j = 1 + \left\lfloor{\log_2\left(\frac{n+1}{2j - 1}\right)}\right\rfloor, ~  \forall j\in\{1,\ldots, k\}.
$$
\end{theorem}
\begin{proof}
A congruence of Lucas (see \S 6.5 in Sagan's textbook 
\cite{Sagan}) states that for a prime $p$ and $k\in\{0,\ldots, m\}$,
\begin{equation}
\label{Lucas}
\binom{m}{k} \equiv \prod\limits_{i\geq 0} \binom{m_i}{k_i}\pmod{p},
\end{equation}
where 
$$
m = \sum\limits_{i\geq 0} m_ip^i ~~~ \mbox{and} ~~~ k = \sum\limits_{i\geq 0} k_ip^i
$$
are the base $p$ expansions of $m$ and $k$ respectively (i.e. $0\leq m_i,k_i < p$).
For prime $p=2$, this congruence implies that for any natural $q\geq 1$, 
\begin{equation}
\label{Lucas2}
(1 + t^m)^{2^q} = \sum\limits_{k = 0}^{2^q} \binom{2^q}{k} t^{mk} 
\equiv 1 + t^{m2^q} \pmod{2}.
\end{equation}
Since, for an arbitrary even $L\leq n+1$ we can find an odd $l\geq 1$ 
such that $L = l2^q$ for some $q\geq 1$, \ref{Lucas2}
implies that in $T\A_{n+1}$, $(1 + t^l)^{2^q} = 1 + t^L$. Therefore, 
the Riordan arrays $a_m = (1 + t^m,\, t)$ with odd $m = 2k - 1$ where 
$k\in\{1,\ldots, \lfloor{(n+1)/2}\rfloor\}$, generate the group $T\A_{n+1}$. 
The minimality of this set of generators follows from the fact that the 
derivative of $1 + t^m\in\F_2[t]$ doesn't have 1 as a root, 
while $1+t^{2k-1}$ has 1 as a root for every $k\geq 1$. 
This proves the statement about the number of factors in the 
direct product decomposition (\ref{DecompA2}).

The kernel of $P_n$ is determined by the single entry in the 
lower left corner of the arrays from $T\A_{n+1}$. It means that 
we can use $a_{n+1} = (1 + t^{n+1},\,t)$ to represent the generator 
of $\ker(P_n)\cong\inte/2\inte$. The identity 
$(1+t^{n+1})^2 \equiv 1\pmod{t^{n+2}}$ 
in $\F_2[t]$ implies $a^2_{n+1} = I$, and the question is if 
$a_{n+1}$ can be written as a product of matrices lifted 
to $T\A_{n+1}$ from $T\A_n$ (i.e. the $(n+2)\times(n+2)$ matrices with 
zero in the lower left corner). Here, by the lifting $s: T\A_n\to T\A_{n+1}$, 
I mean a {\em function} (i.e. not a priori a homomorphism), 
which assigns to the $(n+1)\times (n+1)$ array 
$(p(t),\,t)\in T\A_n$ given by a polynomial $p(t)$ of degree $d\leq n+1$, 
exactly the same pair $(p(t),\,t)$ considered now as an  
$(n+2)\times (n+2)$ array in $T\A_{n+1}$. 
In fact, as will be explained later, the 
lifting $s$ will be a homomorphism if and only if $n$ is even. 

For each odd $m = 2j - 1$ with $j\in\{1,\ldots, \lfloor{(n+1)/2}\rfloor\}$ 
we can find the unique integer $p_m\geq 1$ s.t. 
\begin{equation}
\label{Ineq1}
m2^{(p_m - 1)}\leq n+1 < m2^{p_m}.
\end{equation}
Applying \ref{Lucas2} we obtain the identities 
\begin{equation}
\label{Comp1} 
(1 + t^m)^{2^{(p_m -1)}} = \left\{
\begin{array}{lcl}
1 + t^{n+1}  & \mbox{if} & m2^{(p_m - 1)} = n+1 \\
1 + t^{m2^{(p_m - 1)}} & \mbox{if} & m2^{(p_m - 1)} < n+1. \\
\end{array}
\right.
\end{equation}
The first line in \ref{Comp1} says that when $n$ is odd and $m$ satisfies 
$n + 1 = m2^{(p_m -1)}$ with $p_m\geq 2$, then 
$a_{n+1} = a_m^{2^{(p_m - 1)}}$, that is the array 
$a_{n+1} = (1 + t^{n+1},\,t)$ is not a generator in 
$T\A_{n+1}$. In other words, the short exact sequence \ref{AESZ2}
doesn't split. The identity \ref{Comp1} also implies 
that for every odd $m\in\{1,\ldots, n \}$, the order of $a_m$ in 
$T\A_{n+1}$ will be exactly $2^{p_m}$, since 
$$
m2^{(p_m-1)} \leq n+1 < m2^{p_m}.
$$
Dividing by $m$, and taking $\log_2$ gives
$p_m - 1 \leq \log_2\bigl((n+1)/m\bigr) <  p_m$, i.e.
\begin{equation}
\label{log2}
p_m = p_{2j-1} = 1 + \left\lfloor{\log_2\left(\frac{n + 1}{2j - 1}\right)}\right\rfloor, 
\end{equation}
and the theorem is proven for odd $n$. To prove the case when $n$ is even, 
we show that the splitting function $s:T\A_n\to T\A_{n+1}$ respects all the 
defining relations in the group $T\A_n$, and hence is a homomorphism. 
Now we know that if $n$ is even, 
$T\A_n$ is generated by the arrays $a_{2j-1}, j\in\{1,\ldots, n/2\}$ with the 
corresponding orders $2^{p_{2j-1}}$ given by \ref{log2}. Therefore we need 
to check if $a_{2j-1}^{2^{p_{2j-1}}} = (1,\,t)$ in $T\A_{n+1}$ for all 
$j\in\{1,\ldots, n/2\}$. Using \ref{Lucas2} one more time, we see that in $\F_2[t]$
$$
(1 + t^{2j-1})^{2^{p_{2j-1}}} = 1 + t^{(2j-1){2^{p_{2j-1}}}} \equiv 1 \pmod{t^{n+2}},
$$
since $(2j-1)2^{p_{2j-1}} = m2^{p_m} > n+1$. Thus, the short exact sequence 
\ref{AESZ2} splits, and $T\A_{n+1}\cong T\A_n\times \inte/2\inte$. Using 
these generators with their orders in $T\A_n$, we can write a group presentation 
for the abelian group $T\A_{n+1}$ (skipping all the commutators, since the group 
is commutative) as 
$$
\bigl\langle a_{2j-1}, j\in\{1,\ldots, (n+2)/2\}~|~a^2_{n+1} = 
I, a_{2j-1}^{2^{p_{2j-1}}} = I, 
j\in\{1,\ldots, n/2\} \bigr\rangle.
$$
To prove the formula for $r_k$ with $j= k=  \lfloor{(n+1)/2}\rfloor$, take any 
even $n = 2w\geq 2$, then $2k -1 = 2w-1$. If $w\geq 2$, then 
$$
p_{2k-1} = 1 + \left\lfloor{\log_2\left(\frac{2w+1}{2w - 1}\right)}\right\rfloor = 1,
$$
and if $w = 1$, then we get $p_1 = 2$, as stated in the theorem.
\end{proof}

\noindent Table \ref{tab:INfactors} lists the invariant factors for the 
first few truncated Appell subgroups. 

\begin{table}[h!]
\begin{center}
\begin{tabular}{|l|c|}
\hline
n & Invariant factors of $T\A_{n}(\F_2)$\\
\hline
\hline
1 & \{2\}\\
\hline
2 & \{4\}\\
\hline
3 & \{4,2\}\\
\hline
4 & \{8,2\}\\
\hline
5 & \{8,2,2\}\\
\hline
6 & \{8,4,2\}\\
\hline
7 & \{8,4,2,2\}\\
\hline
8 & \{16,4,2,2\}\\
\hline
9 & \{16,4,2,2,2\}\\
\hline
10 & \{16,4,4,2,2\}\\
\hline
\end{tabular}
\caption{\label{tab:INfactors} Invariant factors of $T\A_{n}(\F_2)$ for $n\in\{1,\ldots,10\}$.}
\end{center}
\end{table}

\begin{corollary}
\label{CorSNAKE}
$\A/\A_1 \cong \inte_2$, and for all $n\geq 1$, 
$$
\A_{2n-1}/\A_{2n+1} \cong \inte_2\times \inte_2.
$$
\end{corollary}
\begin{proof}
Isomorphism $\A/\A_1 \cong \inte_2$ is obvious. 
As follows from our proof of Theorem \ref{DecompA}, the 
kernel of $P_{2n-1}\circ P_{2n}: T\A_{2n + 1} \to T\A_{2n - 1}$ 
is generated by the Riordan arrays $(1 + t^{2n},\,t)$, and 
$(1 + t^{2n + 1},\,t)$. Since each of them has order 2 in $T\A_{2n + 1}$, 
we have $\ker(P_{2n-1}\circ P_{2n})\cong (\inte_2)^2$. Now   
apply the snake lemma to the diagram below.

\begin{equation}
\label{snake}
\begin{tikzcd}
~ & 0 \arrow[d]  & 0 \arrow[d] \arrow[r] & (\inte_2)^2 \arrow[d] & \\
0 \arrow[r] & \A_{2n+1} \arrow[r] \arrow[d, "\pi"] & 
\A \arrow[r] \arrow[d, "\equiv"] & 
T\A_{2n+1} \arrow[d, "P_{2n-1}\circ P_{2n}"] \arrow[r] & 0 \\
0 \arrow[r] & \A_{2n-1} \arrow[r] \arrow[d] & 
\A \arrow[r] \arrow[d] & T\A_{2n-1} \arrow[d] \arrow[r] & 0 \\
~ & \mbox{Coker}(\pi) \arrow[r] & 0  & 0 & ~ \\
\end{tikzcd}
\end{equation}

\end{proof}

Next we prove the main theorem.

\begin{theorem}
\label{LCSRior2}
For all integers $n\geq 2$, 
\begin{equation}
\label{ISOLCSRior2}
\gamma_n(\R) \cong \A_{2n-3} \ltimes \gamma_n(\N),
\end{equation}
and the lower central quotients are
\begin{equation}
\label{ISOLCSRior2LCQ}
\gamma_i(\R)/\gamma_{i+1}(\R) \cong
\left\{ \begin{array}{cl}
(\inte_2)^3\times \inte_4 & \mbox{if $i=1$ } \\
(\inte_2)^4 & \mbox{if $i \geq 2$ is even } \\
(\inte_2)^6 & \mbox{if $i > 1$ is odd. } \\
\end{array}
\right.
\end{equation}
The abelianization of $\R$ has the following group presentation in 
generators and relations 
\begin{equation}
\label{PRESENTR}
 \langle a,\,b_1,\, b_1b_2,\, b_4 ~|~a^2 = (b_1b_2)^2 = b_4^2 = b_1^4 = 1, ~ 
[b_i,\, b_j] =1, ~ [a,\,b_i] = 1 \rangle,
\end{equation}
where $a$ and $b_i$ are the images of the corresponding elements
$$
a_1 = (1 + t,\,t), ~~ e_i = (1,\,t + t^{i+1}).
$$
by the abelianization homomorphism $F_{\R}:\R\to\R^{ab}$.
\end{theorem}

\begin{corollary}
The Riordan group over the field $\F_2$ has width 6.
\end{corollary}

Recall that for all $k\geq 1$ we denote the Riordan array $(1 + \alpha t^k,\,t)$ 
by $a_k(\alpha)$, and the array $(1,\, t + t^{k+1})$ by $e_k$, with 
$a_0 = e_0 = (1,\,t)$. Following the Definition 3 of 
\cite{Camina2} and Definition 12.4.6 of \cite{Leedham}, we will say 
that an element $(1,\,t + t^{m+1} + \cdots)\in \N_m$ has 
{\sl depth $m$}. The identity has depth $\infty$. 
To prove Theorem \ref{LCSRior2} we will need 

\begin{lemma}
\label{CommutA2n}
For all $n\geq 1$ and $i \geq 1$, $[a_n, e^{-1}_i] \in \A$, and 
$$
[a_n, e^{-1}_i] \equiv a_{n + i}[\bar{n}] \pmod{ t^{n + i + 1}},
$$
where $\bar{n}\in\{0,1\}$ and $\bar{n} \equiv n\pmod{2}$.
\end{lemma}
\begin{proof}
It is clear that $[a_n, e^{-1}_i] \in \A$. If we denote the compositional 
inverse of $t + t^{i+1}$ by $\overline{t + t^{i+1}}$, then we have $\pmod{2}$ 
$$
[a_n, e^{-1}_i] = \left(\frac{1}{1 + t^n},\, t + t^{i+1}\right) 
\left(1 + t^n,\, \overline{t + t^{i+1}} \right) 
$$
$$
= \left(\frac{1 + (t + t^{i+1})^n}{1 + t^n},\, t \right) 
= \left(\frac{1 + t^n + t^n
\left(\sum\limits_{s=1}^{n} {n \choose s} t^{si}\right)}{1 + t^n},\, t \right)
$$
$$
= \bigl(1 + t^n(nt^i + {\cal O}(t^{2i}))(1 + t^n + {\cal O}(t^{2n})),\,t\bigr)
$$
$$
\equiv \bigl(1 + nt^{n + i}, \, t\bigr) \pmod{ t^{n+1 + i}} 
$$
as required.
\end{proof}

\begin{lemma}
\label{3rdComm} 
For all integers $n\geq 1$,
$$
\A \cap\gamma_{n+1}(\R) = \A_{2n-1}.
$$
\end{lemma}
\begin{proof} 
We use induction on $n$. Repeating computations from 
Lemma \ref{lemCOM1}, one easily finds that for any $m\geq 2$,
$$
\Bigl[\bigl(1 + t,\,t\bigr),\,\bigl(1,\,t + t^m + t^{m+1}\bigr)^{-1}\Bigr] 
= \bigl( 1 + t^m,\,t\bigl) = a_m,
$$
which implies that $\A_1\leq \A\cap\gamma_{2}(\R)$. 
Considering the projection onto a commutative group 
$\Pi_1:\R \to T\R_1\cong \inte_2$, we see that $a_1\notin \gamma_2(\R)$, 
and hence $\A\cap\gamma_{2}(\R) \leq \A_1$. Assume now 
that the statement holds true for all $n\in\{1,\ldots, L-1\}$, so 
in particular, $a_{2L-1}\in \A_{2L-3} = \A \cap\gamma_{L}(\R)$. 
Applying Lemma \ref{CommutA2n} to the commutator 
$$
\bigl[a_{2L-1},\bigl(1,\, t + t^2\bigr)^{-1} \bigr] \in \bigl[\gamma_L(\R),\,\R\bigr] = 
\gamma_{L+1}(\R),
$$
we deduce that 
$a_{2L} \in \gamma_{L+1}(\R) \pmod {t^{2L + 1}}$. If  
$$
\bigl[a_{2L-1},\bigl(1,\, t + t^2\bigr)^{-1} \bigr] = 
\bigl(1 + t^{2L} + t^{2L+k_1} +\ldots,\,t\bigr),
$$
for some $k_1>0$, take the commutator
$$
\bigl[a_{2L-1},\bigl(1,\, t + t^{2+k_1}\bigr)^{-1} \bigr] = 
\bigl(1 + t^{2L+k_1} +\ldots,\,t\bigr),
$$
and consider the product 
\begin{equation}
\label{ComPr1}
\bigl[a_{2L-1},\bigl(1,\, t + t^2\bigr)^{-1} \bigr] 
\bigl[a_{2L-1},\bigl(1,\, t + t^{2+k_1}\bigr)^{-1} \bigr] 
\end{equation}
$$
\equiv
\bigl(1 + t^{2L} ,\,t\bigr) \pmod{t^{2L+k_1+1}}.
$$
If the product in \ref{ComPr1} equals 
$\bigl(1 + t^{2L} + t^{2L + k_2} +\ldots,\,t\bigr)$ for some $k_2 > k_1$, 
then take the product of three commutators
$$
\prod\limits_{i\in \{0<k_1<k_2\}} \bigl[a_{2L-1},\bigl(1,\, t + t^{2 + i}\bigr)^{-1} \bigr],
$$
and so on. Since the group $\bigl[\gamma_L(\R),\,\R\bigr]$ 
is closed in pro-$p$ topology, we can take the limit of the 
Cauchy sequence made by such products (see Exercise (13) in  
\S 1.6, \cite{Wilson} for the definition), and conclude that 
$a_{2L} = \bigl(1 + t^{2L},\,t\bigr)\in \gamma_{L+1}(\R)$. 
Considering the commutators $\bigl[a_{2L-1},\,(1,\,t+t^m)\bigr]$ for 
$m\geq 3$, one can show similarly that 
$a_n \in \gamma_{L+1}(\R)$ for all $n \geq 2L$, so 
$\A_{2L-1} \leq \A\cap\gamma_{L+1}(\R)$. 

To show that $\A\cap\gamma_{L+1}(\R) \leq \A_{2L-1}$, it will be enough 
to prove that $a_{2L-i} \notin \bigl[\R,\,\gamma_{L}(\R)\bigr]$ 
for all $i\geq 1$. Suppose to the contrary that, for example 
$a_{2L-1} \in \bigl[\R,\,\gamma_{L}(\R)\bigr]$. Then for any integer 
$n > 2L$ there exists $m\in\natu$, such that 
$$
a_{2L-1} \equiv \prod\limits_{i=1}^m [x_i,\,y_i] \pmod{t^n},
$$
where
$$
x_i\in\R, y_i\in \gamma_{L}(\R), i \in \{1,\ldots, m\}.
$$ 
Since $T\A_{2L-1}$ is abelian, there must exist a commutator 
$[x,y] \in \bigl[\R,\,\gamma_{L}(\R)\bigr]$, such that 
$$
\rho_{2L-1}\circ\Pi_{2L-1}(x) \neq 1 ~~ 
\mbox{or} ~~ \rho_{2L-1}\circ\Pi_{2L-1}(y) \neq 1
$$
in the subgroup $\gamma_L(T\N_{2L-1})$ of the truncated Nottingham 
group $T\N_{2L-1}$. But according to Proposition 12.4.30, \cite{Leedham} 
the lowest depth of an element in $\gamma_L(\N)$ equals $2L-1$, that is
$P_{2L-1}(\gamma_L(\N))\cap T\N_{2L-1} = (1,\,t)$. This 
contradiction finishes the proof.
\end{proof}

Proposition 12.4.30 also gives a group 
presentation of the abelianization $\N^{ab} \cong (\inte_2)^2\times \inte_4$ as 
\begin{equation}
\label{PRESENTN}
\N^{ab} \cong \langle b_1,\, b_1b_2,\, b_4 ~|~ (b_1b_2)^2 = b_4^2 = b_1^4 =1, ~ 
[b_i,\, b_j] =1 \rangle,
\end{equation}
where $b_i$ stands for the image of $e_i$ in the quotient group, that is 
$$
b_i := F_{\N} \Bigl( (1,\,t + t^{i+1}) \Bigr) \in \N^{ab}.
$$
Hence, if we take $\D=\F_2$ in the commutative diagram \ref{BCD3ORIG}, 
and use Corollary \ref{CorSNAKE} together with \ref{PRESENTN}, 
we immediately obtain a group presentation for $\R^{ab}$ 
$$
 \langle a,\,b_1,\, b_1b_2,\, b_4 ~|~a^2 = (b_1b_2)^2 = b_4^2 = b_1^4 = 1, ~ 
[b_i,\, b_j] =1, ~ [a,\,b_i] = 1 \rangle,
$$
where $a$ and $b_i$ are the corresponding images by $F_{\R}$, i.e.
$$
a:=  F_{\R} \Bigl( (1 + t,\,t) \Bigr), ~~
b_i := F_{\R} \Bigl( (1,\,t + t^{i+1}) \Bigr).
$$
In particular, we have
$$
\R^{ab} \cong \inte_2\times \N^{ab} \cong (\inte_2)^3\times \inte_4.
$$
That proves the first isomorphism in \ref{ISOLCSRior2LCQ}, 
and \ref{PRESENTR}.

\begin{corollary}
\label{ANORCOR}
For all integers $n\geq 2$, 
$$
\A_{2n-3}\cap\gamma_{n+1}(\R) = \A_{2n-1}.
$$
\end{corollary}
\begin{proof}
Since for any $k\in\{2,\ldots,2n\}$, $\A_{2n-1}\leq \A_{2n-k}$, 
and from Lemma \ref{3rdComm} $\A_{2n-1}\leq \gamma_{n+1}(\R)$, 
the required equality follows from 
$$
\A_{2n-1} \leq \A_{2n-k}\cap\gamma_{n+1}(\R) 
\leq \A \cap\gamma_{n+1}(\R) = \A_{2n-1}.
$$
\end{proof}

Starting with the top exact sequence in the diagram \ref{BCD3ORIG}
\begin{equation}
\label{ESC2}
\begin{tikzcd}
0 \arrow[r] & \A_1 \arrow[r, "\mu^{(2)}"] & 
\gamma_2(\R) \arrow[r, "\rho^{(2)}"]  & 
\gamma_2(\N) \arrow[l, bend left=20, "s^{(2)}"] \arrow[r] & 1,\\
\end{tikzcd}
\end{equation}
using Corollary \ref{ANORCOR}, and repeating the steps we used to 
derive \ref{BCD3ORIG}, one obtains in exactly the same way the following 
commutative diagram.

\begin{equation}
\label{BCD4}
\begin{tikzcd}
& 0 \arrow[d] & 1 \arrow[d] & 1 \arrow[d] & \\
0 \arrow[r] & \A_3 \arrow[r, "\mu^{(3)}"] \arrow[d] & 
\gamma_3(\R) \arrow[r, "\rho^{(3)}"] \arrow[d, "\iota_{\R}"] & 
\gamma_3(\N) \arrow[l, bend left=20, "s^{(3)}"] \arrow[d, "\iota_{\N}"] \arrow[r] & 1 \\
0 \arrow[r] & \A_1 \arrow[r, "\mu^{(2)}"] \arrow[d, "F_{\A}"] & 
\gamma_2(\R) \arrow[r, "\rho^{(2)}"] \arrow[d, "F_{\R}"] & 
\gamma_2(\N) \arrow[l, bend left=20, "s^{(2)}"] \arrow[d, "F_{\N}"] \arrow[r] & 1 \\
0 \arrow[r] & \A_1/\A_3 \arrow[r, "\mu^{}"] \arrow[d] & 
\R_{2/3} \arrow[r,"\rho^{(2/3)}"]  \arrow[d]
& \N_{2/3} \arrow[l, bend left=20, "s^{(2/3)}"] \arrow[r] \arrow[d] & 0 \\
& 0 & 0 & 0 & \\
\end{tikzcd}
\end{equation}

Symbols $\R_{n/(n+1)}:=\gamma_n(\R)/\gamma_{n+1}(\R)$ and 
$\N_{n/(n+1)}:=\gamma_n(\N)/\gamma_{n+1}(\N)$ 
stand for the corresponding lower central quotients of $\R$ and $\N$ 
respectively. Corollary \ref{CorSNAKE} gives 
$\A_1/\A_3 =\langle a_2,a_3\rangle \cong \inte_2\times \inte_2$, and since 
$\N_{2/3}$ is elementary abelian of rank 2 (see again Proposition 12.4.30 in  
\cite{Leedham}) we obtain
$$
\gamma_2(\R)/\gamma_3(\R)\cong (\inte_2)^4.
$$

Repeating this argument we finish the proof of Theorem \ref{LCSRior2} 
by induction. 

\begin{proof}[Proof of Theorem \ref{LCSRior2}]
If $n = 2$, the isomorphism 
$\gamma_n(\R) \cong \A_{2n-3} \ltimes \gamma_n(\N)$ follows from the 
exact sequence \ref{ESC2}. Assume next that for all 
$k\in\{1,\ldots, n \}$, $\gamma_k(\R) \cong \A_{2k-3} \ltimes \gamma_k(\N)$. 
In particular, we have the splitting short exact sequence 

\begin{equation}
\label{ESC3}
\begin{tikzcd}
0 \arrow[r] & \A_{2n-3} \arrow[r, "\mu^{(n)}"] & 
\gamma_n(\R) \arrow[r, "\rho^{(n)}"]  & 
\gamma_n(\N) \arrow[l, bend left=20, "s^{(n)}"] \arrow[r] & 1\\
\end{tikzcd}
\end{equation}

\noindent Using this sequence together with Corollary \ref{ANORCOR}, 
one derives the following commutative diagram, which is 
completely analogous to \ref{BCD4}. 

\begin{equation}
\label{BCD5}
\begin{tikzcd}
& 0 \arrow[d] & 1 \arrow[d] & 1 \arrow[d] & \\
0 \arrow[r] & \A_{2n-1} \arrow[r, "\mu^{(n+1)}"] \arrow[d] & 
\gamma_{n+1}(\R) \arrow[r, "\rho^{(n+1)}"] \arrow[d, "\iota_{\R}"] & 
\gamma_{n+1}(\N) \arrow[l, bend left=20, "s^{(n+1)}"] \arrow[d, "\iota_{\N}"] \arrow[r] & 1 \\
0 \arrow[r] & \A_{2n-3} \arrow[r, "\mu^{(n)}"] \arrow[d, "F_{\A}"] & 
\gamma_n(\R) \arrow[r, "\rho^{(n)}"] \arrow[d, "F_{\R}"] & 
\gamma_n(\N) \arrow[l, bend left=20, "s^{(n)}"] \arrow[d, "F_{\N}"] \arrow[r] & 1 \\
0 \arrow[r] & \A_{2n-3}/\A_{2n-1} \arrow[r, "\mu^{}"] \arrow[d] & 
\R_{n/(n+1)} \arrow[r,"\rho^{(n/(n+1))}"]  \arrow[d]
& \N_{n/(n+1)} \arrow[l, bend left=20, "s^{(n/(n+1))}"] \arrow[r] \arrow[d] & 0 \\
& 0 & 0 & 0 & \\
\end{tikzcd}
\end{equation}

\noindent The top exact sequence in this diagram proves the induction step, 
and hence the isomorphism \ref{ISOLCSRior2}, for all $n\geq 1$.  
Furthermore, since by Proposition 12.4.30, \cite{Leedham} 
we have 
$$
\gamma_i(\N)/\gamma_{i+1}(\N) \cong 
\left\{ \begin{array}{rc}
(\inte_2)^2 & \mbox{if $i$ is even} \\
(\inte_2)^4 & \mbox{if $i$ is odd}, \\
\end{array}
\right.
$$
taking the direct product of each these quotients 
with $\A_{2n-3}/\A_{2n-1} \cong (\inte_2)^2$, proves 
formula \ref{ISOLCSRior2LCQ}, and finishes the proof.
\end{proof}


\subsection{Dihedral groups $D_{2^q}$ as subgroups of $T\R_n$}

It is known that every finite 2-group can be embedded in $\N(\F_2)$, 
and hence in $\R(\F_2)$ (Corollary 12.4.11, \cite{Leedham}). 
A simple argument shows that for any embedding of a finite 
group $G\hookrightarrow \R(\D)$, there is an induced embedding 
$G\hookrightarrow T\R_n(\D)$ for large enough $n$ 
(Proposition 2, \cite{Krylov2}). Therefore, for every 
$q\geq 1$ there exists $n\in\natu$, and an embedding 
$\mu_n:D_{2^q} \hookrightarrow T\R_n(\F_2)$. Here are 
a few examples of such embeddings, for small $q$. 
\begin{itemize}
\item If $q = 0$, $n = 1$: $T\R_1\cong \inte_2\cong D_1$

\item If $q \in\{1,2\}$, $n=2$: For any ring $\D$, $T\R_2(\D)$ is 
isomorphic to the Heisenberg group of upper unitriangular $3\times 3$ 
matrices with coefficients in $\D$, and $T\R_2(\F_2) \cong D_4 > D_2$.

\item If $q =3$, $n=4$: Direct computations show that the group 
$T\R_3$ has 19 nontrivial involutions, and 12 elements of order 4, 
which makes it isomorphic to the central product $D_4 \circ_{Z} D_4$ 
of two copies of $D_4$ over the center. Since $D_8$ has an 
element of order 8, $D_8\nless T\R_3$. However, we can 
embed $D_8$ in $T\R_4$. Indeed, taking 
$r:=(1+t,\,t) \in T\A_4$ as the element of order 8, and 
$s:=(1,\,t+t^2+t^3+t^4)$ as the involution, one easily computes that 
the subgroup $\langle r,\,s\rangle<T\R_4$ generated by $r$ and $s$ 
has 16 distinct elements $\{1,r,r^2,\ldots, r^7,s,sr,sr^2,\ldots,sr^7\}$ and 
hence is isomorphic $D_8$. Here are the finite Riordan array 
representations of these two generators.
$$
r = \begin{pmatrix}
1 & 0 & 0 & 0 & 0\\
1 & 1 & 0 & 0 & 0\\
0 & 1 & 1 & 0 & 0\\
0 & 0 & 1 & 1 & 0\\
0 & 0 & 0 & 1 & 1\\
\end{pmatrix}, ~~ ~~ s = \begin{pmatrix}
1 & 0 & 0 & 0 & 0\\
0 & 1 & 0 & 0 & 0\\
0 & 1 & 1 & 0 & 0\\
0 & 1 & 0 & 1 & 0\\
0 & 1 & 1 & 1 & 1\\
\end{pmatrix}
$$ 

This embedding is hardly a surprise, since $T\A_4$ is the first truncated Appell 
subgroup containing an element of order 8 (recall Table \ref{tab:INfactors}). 
The next theorem extends this example.
\end{itemize}

\begin{theorem}
\label{TheLAST}
For every $n\geq 1$, there exists an embedding $D_{2^{n+1}}\hookrightarrow T\R_{2^n}$. 
\end{theorem}
\begin{proof}
As the {\sl rotation}, take the Riordan array 
$r:=(1+t,\,t)\in T\A_{2^{n}}$, which is an element of order $2^{n+1}$,
due to Theorem \ref{DecompA2}. For the {\sl reflection}, take 
$s:= \bigl(1,\,t + t^2 + \ldots + t^{2^n}\bigr)\in T\N_{2^n}$. Since in 
$\F_2\bigl[[t]\bigr]$ for 
$$
f(t) = t + \ldots + t^{2^n}= t\frac{1+ t^{2^n}}{1+ t}, ~~ 
f\bigl(f(t)\bigr) = t\frac{1+ t^{2^n}}{1+ t}\cdot
\frac{1 + \left(t\frac{1+ t^{2^n}}{1+ t}\right)^{2^n}}{1 + t\frac{1+ t^{2^n}}{1+ t}}
$$
$$
= \frac{t + {\cal O}(t^{2^n+1})}{1 + t^{2^n+1}} 
=\bigl(t + {\cal O}(t^{2^n+1})\bigr)\bigl(1 + {\cal O}(t^{2^n+1})\bigr)
\equiv t \pmod{t^{2^n + 1}}, 
$$
the elements $r$ and $s$ in $T\R_{2^n}$ have orders 
$2^{(n+1)}$ and $2$ respectively. Furthermore, the congruence 
$$
\bigl(1+t,\, t\bigr)\bigl(1,\,t + t^2 + \ldots + t^{2^n}\bigr)\bigl(1+t,\, t\bigr)  
\equiv \bigl(1 + t ,\,t + t^2 + \ldots + t^{2^n}\bigr)\bigl(1 + t,\,t \bigr)
$$
$$
\equiv \bigl(1,\,t + t^2 + \ldots + t^{2^n}\bigr) \pmod{t^{2^n + 1}}
$$
together with a simple argument confirming that all $2^{(n+2)}$ elements
$$
\bigl\{1,r,r^2,\ldots,r^{2^{(n+1)}-1},s,sr,sr^2,\ldots, sr^{2^{(n+1)}-1}\bigr\}
$$ 
are distinct, proves that the subgroup of $T\R_{2^n}$, generated 
by $r$ and $s$, is isomorphic to 
$D_{2^{n+1}} \cong \langle r,s~|~r^{2^{n+1}} = s^2 = 1,~rsr = s\rangle$. 
\end{proof}

\begin{note}
Taking the infinite Riordan arrays 
$$
r = (1 + t,\,t) ~~ \mbox{and} ~~ s = \left(1,\,\frac{t}{1 + t}\right),
$$
one easily checks that they satisfy the relations $s^2 = (1,\,t)$, and 
$rsr = s$. Thus, we obtain an embedding of the infinite dihedral group 
$D_{\infty}$ in the infinite Riordan group $\R(\F_2)$ that induces 
all the embeddings $D_{2^{n+1}}\hookrightarrow T\R_{2^n}$ 
from Theorem \ref{TheLAST} by truncation.
\end{note}

Jointly with Tian-Xiao He, the author found explicit embeddings of an 
arbitrary dihedral group of order $2n$ into the infinite Riordan group 
over the ring of integers modulo $n$ (Theorem 7, \cite{HeKrylov}). 
It would be interesting to see explicit embeddings of the dihedral groups 
$D_{2^n}$ into the infinite Riordan and Nottingham groups over $\F_2$.



\begin{thebibliography}{90}

\bibitem{Babenko} I.K. Babenko, Algebra, geometry, and topology of the substitution 
group of formal power series, {\it Russian Math. Surveys}, 68, (2013), Issue 1, 1 - 68.

\bibitem{BabenkoBogat} I.K. Babenko, S.A. Bogatyi, On the group of substitutions of 
formal power series with integer coefficients, {\it Izv. RAN Ser. Mat.} 72 (2) (2008), 
39 - 64 (in Russian); translation in: {\it Izv. Math.} 72 (2) (2008), 241 - 264.

\bibitem{Barry} P. Barry, {\it Riordan Arrays: A Primer}, Logic Press, 2016.

\bibitem{Bogataya} S.I. Bogataya, S.A. Bogatyy, Series of commutants of the 
Jennings group ${\cal J} (\inte_2)$, {\it Topology and its Applications}, 169 (2014), 
136 - 147.






\bibitem{Camina2} R.D. Camina, The Nottingham group. {\it New horizons in 
pro-$p$ groups}, Progr. Math., Vol 184 (M. du Sautoy, D. Segal, and 
A. Shalev, eds.), Birkh\"auser,  2000, pp 205 - 221.

\bibitem{Cheon1} Gi-Sang Cheon, Nhan-Phu Chung, and Minh-Nhat Phung, 
Algebraic and topological properties of Riordan groups over finite fields, 
https://arxiv.org/abs/2012.10895

\bibitem{Davenport} D.E. Davenport, S.K. Frankson, L.W. Shapiro, L.C. Woodson, 
An Invitation to the Riordan Group, {\it Enumer. Comb. Appl.}, 
ECA 4:3 (2024), Article \# S2S1.

\bibitem{Dixon} J.D. Dixon, M.P.F. du Sautoy, A. Mann, and D. Segal, {\it Analytic 
pro-$p$ Groups}, 2nd ed. Cambrdge studies in advanced math., 61 CUP, 1999


\bibitem{Jennings} S.A. Jennings, Substitution groups of formal power series, 
{\it Canad. J. Math.}, 1954, 6, 325 - 340.




\bibitem{Johnson} D.L. Johnson, The group of formal power series under substitution, 
{\it J. Austral. Math. Soc. (Series A)}, 45 (1988), 296 - 302.

\bibitem{HeKrylov} T.-X. He, N.A. Krylov, On embeddability of Coxeter groups into 
the Riordan group, {\it The Amer. Math. Monthly} to appear in 2026. Preprint 
is available at https://arxiv.org/abs/2405.10470.

\bibitem{Krylov2} T.-X. He, N.A. Krylov, The symmetric groups $S_n, n\geq 4$, and finite
non-abelian simple groups are not embeddable in any Riordan group, Preprint 
is available at https://arxiv.org/abs/2508.03056.

\bibitem{Leedham} C.R. Leedham-Green, S. McKay, {\it The Structure of Groups of 
Prime Power Order}, London Math. Soc. Monogr. vol 27, Oxford University Press, 2002.

\bibitem{Luzon16} A. Luz\'on, D. Merlini, M.A. Mor\'on, L.F. Prieto-Mart{\'i}nez, 
R. Sprugnoli, Some inverse limit approaches to the Riordan group, 
{\it Linear Algebra Appl.} 491 (2016), 239 - 262.



\bibitem{Luzon23} A. Luz\'on, M.A. Mor\'on, L.F. Prieto-Mart{\'i}nez, Commutators and 
commutator subgroups of the Riordan group, {\it Adv. Math.} 
428 (2023), Paper No. 109164, 20 pp.


\bibitem{Sagan} B.E. Sagan, {\it Combinatorics: The Art of Counting}, Graduate Studies 
in Mathematics 210, AMS; 2020.

\bibitem{Shapiro1} L.W. Shapiro, S. Getu, W.-J. Woan, L. Woodson, The Riordan group, 
{\it Discrete Appl. Math.} 34 (1991), 229 - 239.

\bibitem{Shapiro3} L.W. Shapiro, Some Open Questions about Random Walks, 
Involutions, Limiting Distributions, and Generating Functions, {it Adv. Appl. Math.} 
27 (2001), 585 - 596.

\bibitem{Shapiro2} L.W. Shapiro, R. Sprugnoli, P. Barry, G. Cheon, T.-X. He, 
D. Merlini, W. Wang, {\it The Riordan Group and Applications}, Springer, 2022.

\bibitem{Wilson} J.S. Wilson, {\it Profinite groups}, London Math. Soc. Monographs 
New Series 19. Clarendon Press, Oxford, 1998.

\end{thebibliography}
\end{document}